%Paper: alg-geom/9502011
%From: TAN SHENG-LI <tan@ictp.trieste.it>
%Date: Tue, 14 Feb 1995 19:29:40 +0100

%&amstex

% Typeset by AmSTeX
% To appear in J. reine angew. Math.

\input amsppt.sty
\magnification=\magstep1
\hsize15truecm
\vsize24truecm
\font\bigrm=cmr10 scaled1200
\font\Bigrm=cmr10 scaled1728

\def\QED{\hfill{\bf Q.E.D.}}
\def\lra{\longrightarrow}

\def\c_t{\chi_{\text{top}}}

\def\v{\vskip0.2cm}
\def\la{\leftarrow}
\def\o{\over}
\def\wt{\widetilde }

\title\nofrills $ $ \\ \\ \\
\Bigrm Height~inequality~of~algebraic~points \\
on~curves~over~functional~fields
\endtitle
\leftheadtext\nofrills{Tan, Height inequality of algebraic points}
\rightheadtext\nofrills{Tan, Height inequality of algebraic points}
\author{\bigrm Sheng-Li Tan }$^*$ \endauthor

\address
 Department of Mathematics,
 East China Normal University,
 Shanghai 200062,
 P.~R.~of China\endaddress

\thanks$\ast$
The author would like to thank the hospitality and financial support of
Max-Planck-Institut f{\"u}r Mathematik in Bonn during this research.
The manuscript is corrected during a visit to Nice University.
This research is partially supported by the National Natural Science
Foundation of China  and by the
Science Foundation of the University Doctoral Program of CNEC.  \endthanks

\endtopmatter
\hsize15.1truecm
\vsize24truecm
\baselineskip 14truept
\parindent18pt

\vskip0.7cm
\noindent
\centerline{\bf Introduction}
\vskip0.4cm

In this paper, we shall give a linear and effective height inequality for
algebraic
points on curves over functional fields.

\v

Let $f:S\longrightarrow C$ be a fibration of a smooth complex projective
surface $S$ over a curve $C$, and denote by $g$ the genus of a general fiber of
$f$. We assume that $g\geq 2$ and  $S$ is relatively minimal with respect to
$f$,
i.e., $S$ has no $(-1)$-curves contained in a fiber of $f$. Let $k$ be the
functional field of $C$, and $\bar k$ its algebraic closure. For an algebraic
point
$P\in S(\bar k)$, we let $E_P$ be the corresponding horizontal curve on $S$.
The geometric canonical height $h_K(P)$ and the geometric logarithmic
discriminant
$d(P)$ are defined as follows.
    $$
       h_K(P)={K_{S/C}E_P\o [k(P):k]},\hskip0.5cm
       d(P)={2g(\wt E_P)-2 \o [k(P):k]},
    $$
where $\wt E_P$ is the normalization of $E_P$, and $[k(P):k]=FE_P$ is the
degree
of $P$. It is a fundamental problem to give an effective bound of height by the
geometric discriminant. Up to now, many height inequalities have been obtained.

\v
    $$\vbox{
\settabs
\+\hskip-1.5cm Esnault-Viehweg \quad &
$h_{K}(P) \leq (20g - 15)/6 \ d(P) + O(1)$\quad &
 \cr
\+\hskip-1.5cm Szpiro, & $h_{K}(P)\leq
8\cdot3^{3g+1}(g-1)^2(d(P)/3^g+s+1+1/3^{3g})$, &
\cr\smallskip
\+\hskip-1.5cm Vojta,  & $h_{K}(P) \leq (8g-6)/3 \ d(P) + O(1)$,  &
\cr\smallskip
\+\hskip-1.5cm Parshin, & $h_{K}(P) \leq (20g - 15)/6 \ d(P) + O(1)$, &
\cr\smallskip
\+\hskip-1.5cm Esnault-Viehweg, & $h_{K}(P) < 2(2g-1)^2 \ (d(P) + s)$, &
 \cr\smallskip
\+\hskip-1.5cm Vojta, & $h_{K}(P) \leq (2 + \epsilon) \ d(P) + O(1)$, &
\cr\smallskip
\+\hskip-1.5cm Moriwaki, & $h_{K}(P) \leq (2g - 1) \ d(P) + O(1)$, &
\cr
}$$

\v\noindent
where $s$ is the number of singular fibers of $f$.
These inequalities can be found respectively in [Sz], [Vo1], [Pa], [EV],
[Vo2] and [Mo]. It is a problem to get an inequality linear in $g$ with
explicit $O(1)$. (cf. Lang's comments on this problem, [La], p.153).
The purpose of this paper is to give such an inequality.

\v

\proclaim{\indent Theorem~A} Let $f: S\lra C$ be a non-trivial fibration of
genus $g\geq 2$ with $s$ singular fibers, and $P\in S(\bar k)$ an algebraic
point. If $f$ is semistable, then
    $$
       h_K(P)\leq (2g-1)(d(P)+s)-K^2_{S/C},
    $$
and the equality holds only if $f$ is smooth, i.e., $s=0$.

\v
If $f$ is non-semistable, then
    $$
       h_K(P)< (2g-1)(d(P)+3s)-K_{S/C}^2.
    $$
\endproclaim

When we compare it with the canonical class inequality,
the term $3s$ in the second inequality seems natural. Vojta
obtains a canonical class inequality for semistable fibrations:
    $$
       K^2_{S/C}\leq (2g-2)(2g(C)-2+s).
    $$
Furthermore, we have shown that if the equality holds, then $f$ is smooth
(cf. [Ta2], Lemma~3.1). In [Ta1], in a quite natural way, we generalized
Vojta's inequality to the non-semistable case:
    $$
       K^2_{S/C}< (2g-2)(2g(C)-2+3s).
    $$

The first step of the proof of Theorem~A is to obtain the first inequality for
rational points by using Miyaoka-Yau inequality. The idea is motivated by
Xiao's proof of Manin's Theorem (i.e., Mordell conjecture over functional
fields),
(cf. [Xi], Corollary to Theorem~6.2.7).
In the semistable case, the height inequality for algebraic points can be
obtained easily through base changes.
The second step is based on the detailed study of the invariants of semistable
reductions [Ta1].  Kodaira-Parshin's trick plays an important role in this
step.

\v

  I would like to thank Prof. S. Lang for encouraging me to find height
inequalities during my stay at Max-Planck-Institut f\"ur Mathematik in Bonn.
I thank also the referee for his valuable suggestions which make the paper
readable. Due to my carelessness, the referee had to verify by himself
all the results in Sect.~1 for lack of my preprints in which the
proofs are given.

\vskip0.7cm\noindent
\centerline{\bf 1. \ Preliminaries}
\vskip0.4cm

Let $f:S\longrightarrow C$ be a fibration of genus $g\geq 2$,
let $F_1,\cdots,F_s$ be the
singular fibers of $f$, and let $B=\sum_{i=1}^sF_i$. First of all,
we consider the embedded resolution of the singularities of $B_{\text{red}}$.
We denote by $K^2_{S/C}$,
$\chi_f=\deg f_*\omega_{S/C}$ and $e_f=\sum_F(\chi_{\text{top}}(F)-(2-2g))$
the standard relative invariants of $f$.

\v
\proclaim{\indent Definition~1.1} {\rm  The {\it embedded resolution} of
the singularities of $B$ is a sequence
   $$
    (S,B)=  (S_0, B_0)\overset{\sigma_1}\to\la
      (S_1, B_1)\overset{\sigma_2}\to\la
      \cdots
       \overset{\sigma_r}\to\la
       (S_r, B_r)=(S',B')
   $$
satisfying the following conditions.

\v    (i) $\sigma_i$ is the blowing-up of $S_{i-1}$ at a singular point
$(B_{i-1,\text{red}}, p_{i-1})$, which is not an ordinary double point.

\v    (ii) $B_{r, \text{red}}$ has at worst ordinary double points
             as its singularities.

\v    (iii) $B_i$ is the total transformation of $B_{i-1}$.}
\endproclaim

\v
It is well-known that embedded resolution exists uniquely.
We use the notations $m_i$ and $\bar m_i$ to denote respectively
the multiplicities of $(B_{i,\text{red}},p_i)$ and
$(\bar B_{i,\text{red}},p_i)$, where $\bar B_{i,\text{red}}$ is
the strict transform of $B_{\text{red}}$ in $S_i$. Then it is obvious
that

    $$
        \bar m_i\geq m_i-2 .                         \leqno\indent(1)
    $$

Now we let $\pi:\wt C\lra C$ be a base change of degree $d$. Let $S_1$ be the
normalization of $S\times_C\wt C$. We can resolve the singularities of
$S_1$ by using embedded resolution of $B$. It goes as follows.

\v
    $$
         \CD
          S_2@>\eta>>S_1'@>\pi_r>>S'\\
          @\vert@V{\tau}VV@VV{\sigma}V\\
          S_2@>>\rho_2>S_1@>>\rho_1>S
         \endCD
    $$

\v\noindent
where $S_1'$ is the normalization of $S_1\times_SS'$ (hence it is also the
normalization of $S'\times_C\wt C$), and $S_2$ is the minimal resolution of
the singularities of $S_1'$. All of the morphisms are induced naturally.
So $S_2$ is also a resolution of $S_1$. We shall call such a $\rho_2$ the
{\it embedded resolution} of the singularities of $S_1$.

\v
Let $f_2:S_2\lra \wt C$ be the induced fibration, $\wt \rho: S_2\lra \wt S$
the contraction of the $(-1)$-curves contained in the fibers of $f_2$.
Then we have an induced fibration $\wt f: \wt S\lra \wt C$, which is a
relatively minimal fibration determined uniquely by $f$ and $\pi$.
We shall call $\wt f$ the {\it pullback fibration} of $f$ under the
base change $\pi$.

\v
   $$\CD
     \widetilde S@<\widetilde\rho<< S_2 @>\rho_2>> S_1 @>\rho_1>>S \\
      @VV\widetilde fV @VVf_2V @VVf_1V@VVfV \\
       \widetilde C @=\widetilde C@=\widetilde C@>>\pi> C
      \endCD
   $$
Let
   $$
        \Pi_2=\rho_1\circ\rho_2:S_2\longrightarrow S. \leqno\indent(2)
   $$

If $\wt f$ is semistable, then we say that $\pi$ is a {\it semistable
reduction}
of $f$. By using Kodaira-Parshin's construction, we shall construct some
special semistable reductions $\pi$.

\v
Note that if $f$ is semistable,
then $\wt S=S_2$, $\wt \rho$ is the identity endomorphism, and $\Pi_2$ is a
morphism of $\wt S$ to $S$ such that
      $$
     K_{\wt S/\wt C}=\Pi_2^*K_{S/C},\hskip0.3cm K_{\wt S/\wt C}^2
     =dK^2_{S/C}.\leqno\indent(3)
      $$

\proclaim{\indent Lemma~1.2} There exists a semistable reduction
$\pi:\wt C\lra C$ of $f$ such that

\v    {\rm(i)} $\pi$ is ramified uniformly over the $s$ critical points of $f$,
and the ramification index of $\pi$ at any ramified point is exactly $e$.

\v    {\rm (ii)} $e$ is divided by all the multiplicities of the components
of $\sigma^*B$, and it can be arbitrarily large.
\endproclaim

In fact, a base change satisfying the above two conditions must be
a semistable reduction. If $b=g(C)>0$, then the existence follows
from Kodaira-Parshin's construction. Now we consider the case $b=0$.
If $s\geq 3$,
then we can construct a base change totally ramified over the $s$ points.
Then the existence is reduced to the case $b>0$. If $s=2$, the existence
is obvious. Since $f$ is non-trivial, we claim that $s$ is at least $2$.
So the desired semistable reductions always exist.

\v
Indeed, if $b=0$, $s=1$, then $f$ is isotrivial (cf. [Be]). Let $p$
be the critical point. Then $f$ is locally trivial over $\Bbb P^1-p=\Bbb C$.
Because $\Bbb C$ is simply connected, $f$ must be trivial over $\Bbb P^1-p$.
Thus $F\times\Bbb P^1$ is birationally isomorphic to $S$ over $\Bbb P^1$.
By the uniqueness of the relatively minimal model (since $g>0$),
we know that $f$ is trivial.

\v
Let $\Sigma$ be a normal surface, and let $\rho:M\lra \Sigma$ be any
resolution of the singularities. If $\Gamma_1,\cdots,\Gamma_t$ are the
exceptional curves of $\rho$, then the {\it rational canonical divisor}
$K_{\rho}=\sum_{i=1}^t\alpha_i\Gamma_i$ of $\rho$ is defined by the
adjunction formulas
   $$
        K_\rho\Gamma_i+\Gamma_i^2=2p_a(\Gamma_i)-2,
         \hskip0.4cm {\text{for }}\ i=1,\cdots ,t.
   $$
In particular, we have $K_{\rho_2}$ and $K_\eta$.

\v
In Definition~1.1, we denote by $\Cal E_i$ the total inverse image of the
exceptional curve of $\sigma_i$ in $S'$.

\proclaim{\indent Lemma~1.3} Let $\pi$ be the semistable reduction constructed
in
Lemma~1.2. Then we have
    $$
     \wt\rho^*K_{\wt S/\wt C}=\Pi^*_2K_{S/C}-\Pi^*_2\left(
     \sum_{i=1}^s(F_i-F_{i,\text{red}})\right)+K_{\rho_2}-D'' ,
\leqno\indent(4)
    $$
where $D''=K_{S_2/\wt C}-\wt\rho^*K_{\wt S/\wt C}$ is an effective divisor
supported on the exceptional set of $\wt \rho$, and
    $$
      -K_{\rho_2}=\eta^*\pi^*_r\left(\sum_{i=1}^r(m_{i-1}-2)\Cal E_i\right).
\leqno\indent(5)
    $$
\endproclaim
We refer to ([Ta1], \S2.1 and \S5) for the proof of this lemma.  Note that
in this case, $\eta$ is the minimal resolution of rational double points
of type $A_n$, so $K_\eta=0$.

\v
In [Ta1], for each (singular) fiber $F$ of $f$, we associate to it three
nonnegative rational numbers $c_1^2(F)$, $c_2(F)$ and $\chi_F$.
\proclaim{\indent Definition~1.4}{\rm  Let $\pi:\wt C\lra C$ be a base
change of degree $d$ ramified over
$f(F)$ and some non-critical points. If the fibers of $\wt f$ over
$\pi^{-1}(f(F))$ are semistable, then we define}
    $$
       c_1^2(F)=K_{S/C}^2-{1\o d}K_{\wt S/\wt C}^2,\hskip0.2cm
       c_2(F)=e_f-{1\o d}e_{\wt f}, \hskip0.2cm
       \chi_F=\chi_f-{1\o d}\chi_{\wt f}.
    $$
\endproclaim

These invariants are independent of the choice of $\pi$, and can
be computed by embedded resolution of $F$. One of them is zero iff $F$
is semistable. Let
    $$
       I_K(f)=K^2_{S/C}-\sum_{F}c_1^2(F), \hskip0.2cm
       I_\chi(f)=\chi_f-\sum_F\chi_F, \hskip0.2cm
       I_e(f)=e_f-\sum_{F}c_2(F),
    $$
where $F$ runs over the singular fibers of $f$. Then $I_K(f)$,
$I_\chi(f)$ and $I_e(f)$ are nonnegative invariants of $f$, and one of
the first two invariants vanishes if and only if $f$ is {\it isotrivial},
i.e., all the nonsingular fibers are isomorphic. Note that if $f$ is
semistable, then these global invariants are nothing but the standard
relative invariants of $f$.

\proclaim{\indent Lemma~1.5} {\rm([Ta1], Theorem~A$'$)} If $\wt f$ is the
pullback fibration of $f$ under a base change of degree $d$, then we have
    $$
       I_K(\wt f)=dI_K(f),\hskip0.2cm
       I_\chi(\wt f)=dI_\chi(f),\hskip0.2cm
       I_e(\wt f)=dI_e(f).
    $$
\endproclaim

In what follows, we consider the computation of $c_1^2(F)$.
For this, we have to introduce an invariant $c_{-1}(F)$ of $F$.
In fact, if $\pi$ is the semistable reduction in Lemma~1.2,
then $c_{-1}(F)$ can be defined as
    $$
        c_{-1}(F)={1\o \deg \pi}\#\{\, \text{curves over } F \text{ contracted
by }
        \wt\rho \, \}.
    $$
Then we have (cf. (4) or [Ta1], Theorem~3.1)
    $$
      c_1^2(F)=4(g-p_a(F_{\text{red}}))+F_{\text{red}}^2+\sum_{p\in F}\alpha_p
        -c_{-1}(F),
    $$
where $\alpha_p=\sum_{i}(m_{i}-2)^2$, and the $m_i$'s come from the embedded
resolution
of the singular point $(F,p)$. In fact, we have proved that
    $$
       \sum_{p\in F}\alpha_p\leq 2p_a(F_{\text {red}}),
    $$
with equality if and only if $p_a(F_{\text{red}})=0$, i.e., $F$ is a tree
of nonsingular rational curves. (cf. [Ta1], Lemma~3.2). Hence we have

\proclaim{\indent Lemma~1.6} If $F$ is a singular fiber of $f$, then
    $$
       c_1^2(F)+c_{-1}(F)\leq 4g-3,
    $$
and if $p_a(F_{\text{red}})>0$, then
    $$
       c_1^2(F)+c_{-1}(F)\leq 4g-4.
    $$
\endproclaim

Finally, we refer to [Hi] for the details of the following
Miyaoka's inequality.

\proclaim{\indent Lemma~1.7} {\rm ([Mi], Corollary~1.3)} If $S$
is a smooth surface such that the canonical divisor $K_S$ is
nef (numerically effective), and $E_1$, $\cdots$, $E_n$ are
disjoint $ADE$ curves on $S$, then for any (reduced) effective
normal crossing divisor $D$ disjoint to the $E_i$'s, we have
     $$
       \sum_{i=1}^nm(E_i)+3\chi_{\text{top}}(D)\leq 3c_2(S)-(K_S+D)^2,
     $$
where $m(E)$ is defined as follows,
     $$
       \eqalign{ m(A_r)&=3(r+1)-{3\over r+1} ,\cr
                 m(D_r)&=3(r+1)-{3\over 4(r-2)}, \hskip0.2cm \text{ for } r\geq
4,                     \cr
                 m(E_6)&=21-{1\over8},                     \cr
                 m(E_7)&=24-{1\over16},                     \cr
                 m(E_8)&=27-{1\over40}.                     \cr}
     $$
\endproclaim

\vskip0.7cm
\noindent
\centerline{\bf 2. \ The proof of Theorem~A for semistable curves}
\vskip0.4cm

First of all, we give some notations. Let $f:S\lra C$ be a relatively
minimal semistable fibration.
We denote by $f^\#: S^\#\lra C$ the corresponding stable model of $f$,
and by $q$ a singular point of the singular fibers of $f^\#$. Then $q$
is a rational double point of type $A_n$ or a nonsingular point of $S^\#$.
Denote by $E_q$ the inverse image of $q$ in $S$, and let $\mu_q$ be the
Milnor number of $(S^\#,q)$, i.e., the number of $(-2)$-curves in
$E_q$. If $\mu_q=0$, i.e., $q$ is a nonsingular point of $S^\#$, then
$(S^\#,q)$ can be thought of as a ``singular'' point of type $A_0$, thus we
have $m(E_q)=0$ (cf. Lemma~1.7).

\proclaim{\indent Theorem~2.1} If $f:S\lra C$ is non-trivial and semistable,
and $P\in S(\bar k)$ is an algebraic point, then
    $$
         h_K(P)\leq (2g-1)(d(P)+s)-K^2_{S/C},
    $$
and if the equality holds, then $f$ is smooth.
\endproclaim
\v
\demo{\indent Proof} {\it Case} I. $P$ is a $k$ rational point. Let $E$ be
the corresponding section of $f$, and $E^\#$ the image of $E$ in $S^\#$.

\v
If $b=g(C)>0$, then we know that
    $$
        K_S\sim K_{S/C}+(2b-2)F
    $$
is nef. Now we want to use Miyaoka's inequality.
If $q\in E^\#$, then $q$ can not be a nonsingular point of $S^\#$.
Let $E_q^0$ be the $(-2)$-curve in $E_q$ intersecting with $E$. Then
    $$
        E_q-E_q^0=E_{q'}+E_{q''}.
    $$
In this case, we replace $q$ by $q'$ and $q''$. Note that
$m(E_q)=3(\mu_q+1)-3/(\mu_q+1)$, and $\mu_q=\mu_{q'}+\mu_{q''}+1$, ($\mu_{q'}$
and $\mu_{q''}$ may be zero). Hence
    $$
      \varepsilon_q :=  m(E_q)-m(E_{q'})-m(E_{q''})
                   ={3\o \mu_{q'}+1}+{3\o \mu_{q''}+1}-{3\o \mu_q+1}.
    $$
Applying Miyaoka's inequality to $D=E$ and
    $$
        \{ E_q\ \vert \ q\not\in E^\#\}\cup \{ E_{q'}, E_{q''} \ \vert \ q\in
E^\# \},
    $$
we have
    $$
        \sum_q m(E_q)+3\chi_{\text{top}}(E)\leq 3c_2(S)-(K_S+E)^2
                  +\varepsilon  \leqno\indent(6)
    $$
where $\varepsilon=\sum_{ q\in E^\# }\varepsilon_q$.
Note that $f$ is semistable, so $e_f$ is the number of singular
points of the singular fibers of $f$, hence we have
    $
        \sum_q(\mu_q+1)=e_f.
    $
Since $h_K(P)=-E^2$, (6) implies that
    $$
        h_K(P)\leq \sum_q{3\over \mu_q+1}+(2g-1)(2b-2)-K^2_{S/C}+\varepsilon.
                                                 \leqno\indent(7)
    $$
Now we consider the base change $\pi: \wt C\lra C$ of degree $d$
constructed in Lemma~1.2. Let $\wt f: \wt S\lra \wt C$ be the pullback
fibration of $f$ under $\pi$, $\wt P$ the corresponding rational point of
$\wt f$. We use the notation $\wt \cdot$ to denote the corresponding objects of
 $\wt f$. Then the following equalities can be verified easily.
    $$ \eqalign{
        &K_{\wt S/\wt C}^2 =d K^2_{S/C},\hskip0.2cm \wt s={d\over e}s,
       \hskip0.2cm\mu_{\wt q}+1=e(\mu_q+1),
       \hskip0.2cm \wt\varepsilon={d\o e^2}\varepsilon, \cr
        &2g(\wt C)-2  = d(2b-2)+d\left(1-{1\o e}\right)s,
       \hskip0.2cm h_K(\wt P)=dh_K(P).\cr}
    $$
Applying (7) to $\wt f$, we have
    $$
        dh_K(P)\leq {d\o e^2}\sum_q{3\o \mu_q+1}+(2g-1)\left((2b-2)d+
                     d\left(1-{1\o e}\right)s\right)
                    -dK^2_{S/C}   +{d\o e^2}\varepsilon,
    $$
i.e.,
    $$
        h_K(P)-(2g-1)\left(d(P)+s\right)+K^2_{S/C}\leq
                      -{(2g-1)s\o e}+{1\o e^2}\left(\sum_q{3\o \mu_q+1}
                     +\varepsilon\right).
    $$
Let $e$ be large enough we can see that the left hand side $\leq 0$,
or $<0$ if $s> 0$.

\v
Now we consider the case $b=0$. Since $f$ is non-trivial,  we have
$s\geq 5$ [Ta2]. Then we consider
also the base change given in Lemma~1.2. Since $g(\wt C)>0$, the height
inequality for $\wt P$ holds, which implies the inequality for $P$.

\v
{\it Case} II. $P$ is an algebraic point of degree $d_P$. Let $E_P$ be
the corresponding reduced and irreducible horizontal curve on $S$, $\wt C$
the normalization of $E_P$,
and $\pi: \wt C\lra C$ the morphism induced by $f$. Let
$\wt f: \wt S\lra \wt C$ be the pullback of $f$ under $\pi$.
Since $f$ is semistable, from (2) and (3) we know that $\Pi_2$ is a morphism
of $\wt S$ to $S$ such that
    $$
       K_{\wt S/\wt C}=\Pi_2^*(K_{S/C}), \hskip0.3cm
       K^2_{\wt S/\wt C}=d_PK^2_{S/C}.
    $$
By the construction of $\wt f:\wt S\lra \wt C$, there is a section
$\wt E$ of $\wt f$ such that the induced map of $\wt E$ to $E_P$ is a
birational morphism, hence ${\Pi_2}_*(\wt E)=E_P$. By projection formula
we have
    $$\eqalign{
       h_K(P)&={1\o d_P}E_PK_{S/C}   \cr
             &={1\o d_P}\wt E\cdot \Pi^*_2(K_{S/C})   \cr
             &={1\o d_P}\wt E K_{\wt S/\wt C}   \cr
             &\leq (2g-1)\left({2g(\wt C)-2\o d_P} + {\wt s\o d_P}\right)
                             -{1\o d_P}K^2_{\wt S/\wt C} \cr
             &\leq (2g-1)(d(P)+s)-K^2_{S/C}.}
     $$
If $s> 0$, then the inequality holds strictly.
\QED
\enddemo

%\newpage
\vskip0.7cm
\noindent
\centerline{\bf 3. \ The proof of Theorem~A for non-semistable curves}
\vskip0.4cm

   Let $f:S\lra C$ be a non-semistable fibration with $s$ singular fibers
$F_1,\cdots, F_s$. Let $P$ be an algebraic point of degree $d_P$, and $E_P$
the corresponding horizontal curve on $S$. We shall prove in this section
that
     $$
        h_K(P)< (2g-1)(d(P)+3s)-K^2_{S/C}.              \leqno\indent(8)
     $$
We let $\pi: \wt C\lra C$ be the semistable reduction of $f$ constructed
in Lemma~1.2. For convenience, we recall the construction of the pullback
fibration $\wt f$ under $\pi$ given in Sect.~1.

   $$
     \CD
               @.   S_2@>\eta>>S_1'@>\pi_r>>S'\\
               @.  @| @VV{\tau}V@VV{\sigma}V\\
      \widetilde S@<\widetilde\rho<< S_2 @>\rho_2>> S_1 @>\rho_1>>S \\
      @VV\widetilde fV @VVf_2V @VVf_1V@VVfV \\
       \widetilde C @=\widetilde C@=\widetilde C@>>\pi> C
      \endCD
   $$
Note that $\Pi_2$ is defined as $\rho_1\circ\rho_2$.

\v
Let $C_P$ be the normalization of $E_P$,  $\pi_P:C_P\lra C$ the morphism
induced by $f$, and $f_P:S_P\lra C_P$ the pullback fibration of $f$ under
$\pi_P$. By the construction of $f_P$, there is a section of $f_P$ which
maps birationally onto $E_P$.

\v
Now considering the normalization of one component of the fiber product
$\wt C\times_C C_P$, we obtain a curve $\hat C$ such that the following
diagram commutes.
     $$     \CD
            \hat C@>\psi>>C_P\\
            @V\phi VV@VV\pi_PV\\
            \wt C@>>\pi>C
            \endCD
     $$
Among all the components of $\wt C\times_CC_P$, we can choose
$\hat C$ such that the induced morphism $\phi$ is of the least degree.

\v
Let $\hat f:\hat S\lra \hat C$ be the pullback fibration of $\wt f$ under
$\phi$. By the uniqueness of the relatively minimal model (since $g>0$)
and the universal property of fiber product, we know that $\hat f$ is
nothing but the pullback of $f_P$ under $\psi$.
     $$    \CD
           \hat S\hskip0.22cm{\hskip0.525cm\over}
           {\hskip0.525cm\over}\hskip-1.4cm@.
           {\hskip0.525cm\over}\hskip-0.15cm@>>>S_P\\
           @V\Phi_2VV@.@VVV\\
           \wt S@<<\wt\rho<S_2@>>\Pi_2>S
           \endCD
     $$
Hence $\hat f$ has a section $\hat E$ induced by the section of $f_P$
mentioned above.
Now we know that the induced rational map of $\hat E$ to $E_P$ is of degree
$\deg\psi$. Since $\wt f$ is semistable, the
induced map $\Phi_2: \hat S \lra \wt S$ is a morphism. Denote by $\wt E$
the image of $\hat E$ in $\wt S$, and let $E_2$ be the strict transform
of $\wt E$ in $S_2$.

\v
We claim that the restriction map $\Phi_2\vert_{\hat E}: \hat E\lra \wt E$ is
birational.
     $$   \CD
          \hat E@>>>\hat S@>\hat f>>\hat C\\
                @V\Phi_2\vert_{\hat E}VV@V\Phi_2VV@VV\phi V\\
          \wt E@>>>\wt S@>>\wt f>\wt C
          \endCD
     $$

\v
Indeed, if $\Phi_2\vert_{\hat E}$ is not birational, then
the degree of the morphism $\wt f|_{\wt E}: \wt E\lra \wt C$ is less than that
of $\phi$. On the other hand, we can see that the map of $\wt E$ to $E_P$
can be lifted to $C_P$, hence $\wt f|_{\wt E}$ is factorized through
 $\wt C\times_C C_P$.
Hence there is a component of $\wt C\times_C C_P$ such that degree of the
induced
morphism to $\wt C$ is less than that of $\phi$, which contradicts the choice
of $\hat C$. This proves the claim.

\v
Now we have
   $$\eqalign{
    {\Pi_2}_*E_2&=\deg\psi\ E_P, \hskip0.3cm {\Phi_2}_*\hat E=\wt E,\cr
    \wt \rho_*E_2&=\wt E,\hskip0.3cm K_{\hat S/\hat C}
    =\Phi_2^*K_{\wt S/\wt C}.}
                                                   \leqno\indent(9)
   $$
{}From Lemma~1.3,
     $$
         \wt \rho^*K_{\wt S/\wt C}=\Pi^*_2K_{S/C}-D_\pi,
     $$
where
     $$
        D_\pi =\Pi_2^*\left(\sum_{i=1}^s(F_i-F_{i,\text{red}})\right)
               -K_{\rho_2}+D''
     $$
is an effective divisor. By (9) and projection formula we have
     $$\eqalign{
        h_K(\hat P)&=K_{\hat S/\hat C}\hat E=\Phi^*_2K_{\wt S/\wt C}\hat E\cr
                   &=K_{\wt S/\wt C}\wt E=\wt\rho^*K_{\wt S/\wt C}E_2\cr
                   &=(\Pi_2^*K_{S/C}-D_\pi)E_2\cr
                   &=\deg\psi \ K_{S/C}E_P-D_\pi E_2\cr
                   &=d_P\deg\psi\ h_K(P)-D_\pi E_2,}
     $$
thus
     $$
        h_K(P)={1\o d_P\deg\psi}h_K(\hat P)+{1\o d_P\deg\psi}D_\pi E_2 .
         \leqno\indent(10)
     $$
Note that
     $$
        d\deg\phi=d_P\deg \psi .                      \leqno\indent(11)
     $$
where $d$ is the degree of $\pi$.

\proclaim{\indent Lemma~3.1}
     $$
        {1\o d_P\deg\psi}h_K(\hat P)\leq (2g-1)(d(P)+s)-I_K(f).
     $$
\endproclaim
\demo{\indent Proof} First of all, we consider the case when
$\hat C=\Bbb P^1$ and $\hat f$ is trivial, i.e., $f$ is isotrivial,
so $I_K(f)=0$. It is easy to see that $h_K(\hat P)=0$.
Since $f$ is non-trivial, $s\geq 2$, and the desired inequality
holds.

\v
In what follows,
we assume that $\hat f$ is non-trivial if $\hat C=\Bbb P^1$.

\v
Since $\hat f$ is semistable, by Theorem~2.1, we have
     $$
       h_K(\hat P)\leq (2g-1)\left({2g(\hat C)-2 }+
          {\hat s }\right)-K^2_{\hat S/\hat C},
     \leqno\indent(12)
     $$
where $\hat s$ is the number of singular fibers of $\hat f$.
(Note that if $\hat f$ is trivial, then we have a stronger
inequality $h_K(\hat P) \leq 2g(\hat C)-2$, because $g(\hat C)>0$).
It is obvious that
     $$
        \hat s\leq {ds\o e}\deg\phi ={s\o e}d_P\deg\psi.   \leqno\indent(13)
     $$
{}From Lemma~1.5, we have
     $$
       K^2_{\hat S/\hat C}={d_P\deg\psi }\ I_K(f).       \leqno\indent(14)
     $$
By Hurwitz formula,
     $$
        2g(\hat C)-2=\deg\psi \ (2g(C_P)-2)+r_\psi.
     $$
Note that the ramification index of $\pi$ at any ramified point is $e$,
by the construction of $\psi$  we can see that the index of $\psi$ at any
ramified point is at most $e$. Let $x$ be a branch point of $\psi$. Then
we know that the contribution of $\psi^{-1}(x)$ to $r_\psi$ is at most
$(1-1/e)\deg\psi$.
On the other hand, $\psi$ has at most $d_Ps$ branch points, thus
     $$
        {r_\psi}\leq \left(1-{1\o e}\right)sd_P\deg\psi,
     $$
it implies that
     $$
       {2g(\hat C)-2 }\leq \left(d(P)+ \left(1-{1\o e}\right)s\right)
                           { d_P\deg\psi }.
                                                       \leqno\indent(15)
     $$
Combining (12)--(15), we have
     $$
       h_K(\hat P)\leq \left((2g-1)(d(P)+s)-I_K(f)\right){ d_P\deg\psi }.
     $$
This completes the proof.
\QED
\enddemo

\v

Now we shall find the upper bound of $D_\pi E_2/{ d_P\deg\psi }$. So
we only need to consider the first diagram at the beginning of this section.
{}From (9),
    $$
{\Pi_2}_*E_2 =\deg\psi\ E_P,\hskip0.4cm E_P(F_i-F_{i,\text{red}})<d_P,
    $$
by projection formula we have
\proclaim{\indent Lemma~3.2}
     $$
       {1\o d_P\deg\psi }\Pi^*_2
       \left(\sum_{i=1}^s(F_i-F_{i,\text{red}})\right) E_2< s.
     $$
\endproclaim

\proclaim{\indent Lemma~3.3}
     $$
       { -K_{\rho_2}E_2\over d_P\deg\psi}
       \leq s-\#\{ F_i\ \vert \ p_a(F_{i, \text{red}})=0\}.
     $$
\endproclaim
\demo{\indent Proof} Note first that $p_a(F_{i,\text{red}})=0$ implies that
$F_i$ is a tree of non-singular rational curves, so $F_i$
has no effect on $-K_{\rho_2}$ and $-K_{\rho_2}E_2$.
For simplicity, we assume that $p_a(F_{i,\text{red}})\neq 0$ for all $i$.

\v
We let
     $$
        \sigma^*E_P=\bar E_P + \sum_{i=1}^ra_{i-1}\Cal E_i,
     $$
where $\bar E_P$ is the strict transform of $E_P$ and $a_i\geq 0$ is the
multiplicity of the strict transform of $E_P$ at $p_i$. We have known that
$(\pi_r\circ\eta)_*E_2=\deg\psi\ \bar E_P$, and
    $$
       -K_{\rho_2}=(\pi_r\circ\eta)^*\left(\sum_{i=1}^r(m_{i-1}-2)\Cal
E_i\right),
    $$
hence
     $$\eqalign{
       -K_{\rho_2}E_2&=(\pi_r\circ\eta)^*\left(\sum_{i=1}^r(m_{i-1}-2)
                          \Cal E_i\right)E_2\cr
                     &=\deg\psi\ \sum_{i=1}^r(m_{i-1}-2)\Cal E_i \bar E_P\cr
                     &=\deg\psi\ \sum_{i=1}^r(m_{i-1}-2)a_{i-1}.\cr}
                       \leqno\indent(16)
     $$
On the other hand,
    $$
       \sigma^*\left(\sum_{i=1}^sF_{i,\text{red}}\right)=\sum_{i=1}^s
       \bar F_{i,\text{red}}+
                   \sum_{i=1}^r\bar m_{i-1}\Cal E_i,
    $$
where $\bar F_{i,\text{red}}$ is the strict transform of $F_{i,\text{red}}$,
and $\bar m_{i-1}$ is
the multiplicity of the strict transform of $\sum_{i=1}^sF_{i,\text{red}}$ at
$p_{i-1}$.
{}From $\sum_{i=1}^s\bar F_{i,\text{red}}\bar E_P\geq 0$, we have
    $$
      \sum_{i=1}^ra_{i-1}\bar m_{i-1}\leq
      \sum_{i=1}^sF_{i,\text{red}}E_P\leq sd_P. \leqno\indent(17)
    $$
Combining (16) and (17) with (1), we have
    $$
       -K_{\rho_2}E_2\leq sd_P\deg\psi .
    $$
This completes the proof of the lemma.
\QED
\enddemo

\proclaim{\indent Lemma~3.4}
     $$
       {D''E_2\o d_P\deg\psi }\leq\sum_{i=1}^sc_{-1}(F_i).
     $$
\endproclaim
\demo{\indent Proof} We let
$\wt\rho=\wt\rho_k\circ\wt\rho_{k-1}\circ\cdots\circ\wt\rho_1$
 be the decomposition of $\wt\rho$ into $k$ blowing-downs of $(-1)$-curves. If
we denote by
$\Delta_i\subset S_2$ the total transform of the exceptional curve of
$\wt \rho_i$, then
   $$
       D''=K_{S_2/\wt C}-\wt \rho^*K_{\wt S/\wt C}=\sum_{i=1}^k\Delta_i.
   $$
It is easy to see that $E_2\Delta_i$ is the multiplicity of
$\wt\rho_i\circ\cdots\circ\wt\rho_1(E_2)$ at
$\wt\rho_i\circ\cdots\circ\wt\rho_1(\Delta_i)$, so
   $$
      E_2\Delta_i\leq E_2F ={\deg\psi\o d}d_P,
   $$
where $F$ is a fiber of $f_2$.
 On the other hand, the number of curves
contracted by $\wt \rho$ is $k=d\sum_{i=1}^sc_{-1}(F_i)$. Hence we have
the desired inequality.
\QED
\enddemo

\vskip0.4cm
{\it Proof of (8).} \ From (10) and the above lemmas, we have
     $$\eqalign{
       h_K(P)  <
&(2g-1)(d(P)+s)-K^2_{S/C}+\sum_{i=1}^s\left(c_1^2(F_i)+c_{-1}(F_i)
                   \right)\cr
                    &-\#\{F_i \ \vert \ p_a(F_{i,\text{red}})=0\}+2s.}
     $$
By Lemma~1.6,
     $$
       \sum_{i=1}^s\left(c_1^2(F_i)+c_{-1}(F_i)\right)\leq
       (4g-4)s+\#\left\{ F_i\ \vert \ p_a(F_{i,\text{red}})=0 \right\}.
     $$
Hence we have
     $$
       h_K(P)< (2g-1)(d(P)+3s)-K^2_{S/C}.
     $$
\QED

\enddocument